\documentclass[11pt]{article}
\usepackage{amsfonts}
\usepackage[mathscr]{eucal}
\usepackage{mathrsfs}
\usepackage{amsmath}
\usepackage{amssymb}

\newtheorem{theorem}{Theorem}[section]
\newtheorem{lemma}[theorem]{Lemma}
\newtheorem{proposition}[theorem]{Proposition}
\newtheorem{definition}[theorem]{Definition}
\newtheorem{remark}[theorem]{Remark}
\newtheorem{corollary}[theorem]{Corollary}

\newtheorem{assumption}[theorem]{Assumption}
\newcommand{\De}{\Delta}

\newcommand{\Om}{{\Omega}}

\newcommand{\e}{\epsilon}
\newcommand{\ga}{{\gamma}}
\newcommand{\eps}{\varepsilon}

\newcommand{\si}{{\sigma}}

\newcommand{\Ac}{\mathfrak{A}}

\newcommand{\sB}{\mathscr{B}}
\newcommand{\sC}{\mathscr{C}}

\newcommand{\N}{{\mathbb N}}

\newcommand{\R}{{\mathbb R}}

\newcommand{\cB}{{\cal B}}

\newcommand{\cO}{{\cal O}}
\newcommand{\cF}{{\cal F}}

\newcommand{\cE}{{\cal E}}
\newcommand{\cH}{{\cal H}}

\newcommand{\cM}{{\cal M}}
\newcommand{\cN}{{\cal N}}

\newcommand{\cL}{{\cal L}}

\newcommand{\g}{{\nabla}}

\newcommand{\pd}{\partial}

\newcommand{\Lto}{\widehat L_2(\Om)}
\newcommand{\Hto}{\widehat H^2_{0}(\Om)}
\newcommand{\hch}{\widehat{\cH}}
\newcommand{\wH}{\widehat{H}}

\newcommand{\dist}{{\operatorname{dist}}}

\newcommand{\di}{{\rm div\, }}
\newcommand{\wrt}{{with respect to }}

\newenvironment{declaration}[1]{\trivlist
\item[\hskip \labelsep{\bf #1 }]\ignorespaces}{\endtrivlist}
\newenvironment{proofof}[1]{\begin{declaration}{#1}}{\hfill
$\square$ \end{declaration}}
\newenvironment{proof}{\begin{proofof}{Proof.}}{\end{proofof}}

\begin{document}
\title{On interaction of an elastic wall \\ with a Poiseuille type flow}
\author{Igor Chueshov\footnote{\small e-mail:
 chueshov@univer.kharkov.ua}
 ~~ and ~~
  Iryna Ryzhkova\footnote{\small e-mail:
 iryonok@gmail.com} \\
 \\Department of Mechanics and Mathematics, \\
 Kharkov National University, \\ Kharkov, 61077,  Ukraine\\  }

\maketitle

\centerline{\it Dedicated to Academician Anatoly M. Samoilenko}
\centerline{\it on the occasion
     of his $75$th birthday}
\begin{abstract}
We study dynamics of a  coupled system consisting of  the 3D
Na\-vier--Stokes equations which is linearized near a certain
Poiseuille type flow in an (unbounded) domain  and a classical
(possibly nonlinear) elastic plate equation for transversal displacement
on a flexible flat  part  of the boundary.
We first show that this problem generates an evolution semigroup $S_t$ on an appropriate phase space.
Then under some conditions concerning the underlying (Poiseuille type) flow
we prove the existence of a compact finite-dimensional global attractor for this semigroup
and also show that  $S_t $ is an exponentially stable $C_0$-semigroup
of linear operators  in the fully  linear case. Since
we do not assume any kind of mechanical damping in the plate component,
this means that dissipation of the energy in the fluid flow
due to viscosity is sufficient to stabilize the system.

\end{abstract}

\section{Introduction}
Let $\cO\subset \R^3$ be a (possibly unbounded) domain  with a sufficiently smooth
 boundary $\partial\cO$.
 We assume that
$\partial\cO=\overline{\Omega}\cup \overline{S}$,
 where $\Om\cap S=\emptyset$,
$$
\Om\subset\{ x=(x_1;x_2;0)\, :\,x'\equiv
(x_1;x_2)\in\R^2\}~~
$$
is bounded in $\R^2$ and has  the smooth contour $\Gamma=\partial\Om$.
We refer to Assumption~\ref{as:dom} below for further hypotheses concerning
the domain $\cO$.

 \par
 Let $a_0(x)=(a_0^1(x);a_0^2(x);a_0^3(x))$ be a smooth bounded field defined on $\overline{\cO}$ such that
 $\di a_0=0$,
$(n,a_0)=0$ on $\pd\cO$ ($n$ is the exterior normal to $\partial\cO$,
$n=(0;0;1)$ on $\Om$)
and $A=A(x)$ be a  bounded  measurable $3\times 3$ matrix,  $x\in \overline{\cO}$.
We introduce   a linear first order operator $L_0$ of the form
\begin{equation}\label{oper-L0}
L_0v= (a_0, \g)v +A v
\end{equation}
 and consider the following {\em linear} Navier--Stokes equations in $\cO$
for the fluid velocity field $v=v(x,t)=(v^1(x,t);v^2(x,t);v^3(x,t))$
and for the pressure $p(x,t)$:
\begin{equation}\label{fl.1}
   v_t-\nu\Delta v+L_0v+\nabla p=G_f(t)\quad {\rm in\quad} \cO
   \times(0,+\infty),
\end{equation}
   \begin{equation}\label{fl.2}
   \di v=0 \quad {\rm in}\quad \cO
   \times(0,+\infty),
  \end{equation}
where $\nu>0$ is the dynamical viscosity, $G_f(t)$ is a volume force
(which may depend on $t$).
We supplement (\ref{fl.1}) and (\ref{fl.2}) with  the (non-slip)  boundary
   conditions imposed  on the velocity field $v=v(x,t)$:
\begin{equation}\label{fl.4}
v=0 ~~ {\rm on}~S;
\quad
v\equiv(v^1;v^2;v^3)=(0;0;u_t) ~~{\rm on} ~ \Om.
\end{equation}
Here, as in \cite{ChuRyz2011}, $u=u(x,t)$ is the transversal displacement
of the plate occupying $\Om$ and satisfying
the following  equation
\begin{equation}\label{pl_eq}
u_{tt} + \De^2 u + \cF(u)=G_{pl}(t) +p|_\Om
~~{\rm in}~~ \Omega \times (0, \infty),
\end{equation}
where $\cF(u)$ is  a nonlinear (feedback) force
(see Assumption~\ref{A:force} below), $p$ is the pressure
from \eqref{fl.1}, $G_{pl}(t)$ is a given external (non-autonomous) load.
We refer to \cite{ChuRyz2011} and to the references therein
 for some discussion of this plate model and
for an explanation  of the structure of the force exerted by the fluid on the plate.
\par
We impose clamped boundary conditions on the plate deflection
\begin{equation}
u|_{\pd\Om}=\left.\frac{\pd u}{\pd n} \right|_{\pd\Om}=0 \label{plBC}
\end{equation}
and supply \eqref{fl.1}--\eqref{plBC} with initial data of the form
\begin{equation}
 v(0)=v_0,\quad u(0)=u_0, \quad u_t(0)=u_1. \label{IC}
\end{equation}
If we assume that the velocity field $v$ decays sufficiently fast as $|x|\to+\infty$ and $x\in\overline{\cO}$,
then  \eqref{fl.2} and \eqref{fl.4} imply the following
compatibility condition
\begin{equation}\label{Com-con}
\int_\Om u_t(x',t) dx'=0 \quad \mbox{for all}~~ t\ge 0,
\end{equation}
which  can be interpreted as preservation of the volume of the fluid.
\par
Below (see Definitions~\ref{lin_de:solution} and \ref{de:solution})
we define  a solution to
\eqref{fl.1}--\eqref{Com-con} as a pair $(v;u)$ satisfying
some variational type relation. If the pair $(v;u)$ is already
determined, then (at least formally) we can find $\g p$ in $\cO$
and the trace of $p$ on $\Om$ from  \eqref{fl.1} and \eqref{pl_eq}.
Thus the pressure $p$ is uniquely defined by $(v;u)$.
\par
The main example which we have in mind is the {\em Poiseuille} flow
(see, e.g., \cite{chorin-marsden} for some details).
In this case we deal with the domain
  \begin{equation}\label{poisl-1}
\cO=\{(x_1;x_2;x_3)\, : (x_2;x_3)\in \sB\subset\R^2,~~x_1\in \R\},
  \end{equation}
  where $\sB$ is a domain in $\R^2$,
and  the Poiseuille velocity field has the form $a_0=(a(x_2;x_3);0;0)$, where $a(x_2;x_3)$ solves
the elliptic problem
\begin{equation}\label{poisl-1a}
\nu\Delta a=-k~\mbox{ in }~\sB,~~~ a=0~\mbox{ on }~\pd \sB,
\end{equation}
where $k$ is a positive parameter.
Linearization  of the nonlinear Navier-Stokes equations around  the flow $a_0$ gives us the model with
\begin{equation}\label{poisl-2}
L_0v = (a_0, \g)v + (v, \g) a_0.
\end{equation}
There are two important special cases of the choice of $\sB$ in \eqref{poisl-1}: (i) $\sB$ is a bounded
 domain in $\R^2$ (the Poiseuille flow in a cylindrical tube) and (ii)
a flow between two parallel planes. In the latter case
\begin{equation}\label{lyier}
\sB=\{(x_2;x_3): x_2\in\R,~x_3\in(-h,0)\},~~a(x_2;x_3)=-\frac{k x_3}{2\nu}(h+x_3).
\end{equation}
\par
Another possibility included in the framework above is the {\em Oseen} modification of \eqref{fl.1}
(see, e.g., \cite{chorin-marsden}).
In this case $L_0v=U\pd_{x_1}v$, where $U$ is the parameter
which has the sense of the speed of the unperturbed flow moving along the $x_1$-axis.
This corresponds to the case when $a_0=(U;0;0)$ and $A(x)\equiv 0$  in \eqref{oper-L0}.
We can also consider the   situation when  $a_0\equiv 0$ and $A(x)\not\equiv 0$  in \eqref{oper-L0}.
In this case we note that if $A(x)$ is a symmetric strictly positive matrix (e.g.,  $A(x)=\si I$, $\si>0$),
then  $L_0v=A(x)v$  can be interpreted as a drag/friction term which models the  resistance
offered by the fluid against its flow (see, e.g., \cite{massey}  for some discussion).

Thus, our general model includes the case of
interaction  of the Poiseuille (or Oseen type) flow (with a possible drag/friction) in the domain $\cO$
 bounded by
the (solid) wall $S$ and a horizontal boundary $\Om$
on which a thin (nonlinear) elastic plate is placed.
The motion of the fluid is described
by the  3D Navier--Stokes equations {\em linearized} around the
Poiseuille (or Oseen) flow $a_0(x)$.
To describe  deformations of the plate
 we consider  a generalized plate model  which accounts only for
transversal displacements and
 covers a general large deflection Karman type model  and
 can be also applied to nonlinear Berger and Kirchhoff plates
(see the discussion in \cite{ChuRyz2011} and also in Section~\ref{sec:nonlin}).
Since we deal with {\em linearized} fluid equations the  interaction
  model considered  assumes that  large deflections of the plate produce
 small effect on the corresponding  underlying  flow.
\par
We note that the mathematical studies of the problem
of fluid--structure interaction
in the case of viscous fluids and elastic plates/bodies
have a long history.
We refer to \cite{CDEG05,Chu_2010,ChuRyz2011,Gra2008,Ggob-aa09,Ggob-mmas09,GPP2008,Kop98,OP1999}
and the references therein for the case of plates/membranes.
The case of moving elastic bodies \cite{CS06}
and  the case of elastic bodies with the fixed interface \cite{avalos-amo07,aval-tri09,BGLT07,DGHL03}
were studied;
see also  the literature cited in these references.
We also mention the recent short survey \cite{berlin11}  and the paper  \cite{cr-full-karman} which
deals with dynamical issues for a model taking into account both transversal
and longitudinal deformations. All these sources deals with the case of bounded reservoirs $\cO$.
\smallskip \par
In this paper our main point of interest is well-posedness and
long-time dynamics of
solutions to the coupled problem in
\eqref{fl.1}--\eqref{IC}
for the velocity $v$ and the displacement $u$ in the case of {\em unbounded}
domains $\cO$.
\par
In our argument we use the ideas and methods developed in
our previous paper \cite{ChuRyz2011}.
Our  main difficulties in comparison with  \cite{ChuRyz2011}
is related to the facts that (i) we  deal with the (possibly) {\em unbounded} domain $\cO$
(hence, we loose some compactness properties of
the fluid velocity variable and cannot use eigenfunctions of the Stokes operator)
 and (ii) the  fluid equation \eqref{fl.1} is perturbed by nonconservative and nondissipative term
(hence, we  can loose the energy monotonicity and need some additional argument
for non-monotone parts).
To overcome these difficulties we are enforced to
 use a general basis in the fluid component
and a specially constructed extension operator  $Ext$
of functions on  $\Om$  into solenoidal functions on $\cO$.
\par
The paper is organized as follows.
In Section~\ref{sec:pre} we introduce Sobolev type spaces we need
 and provide with some results concerning the extension operator
$Ext$.
In Section~\ref{sec:lin} we prove Theorem~\ref{lin_WP} on  well-posedness in the case of
linear model and study stability properties of solutions in Theorem~\ref{th:stab}.
Section~\ref{sec:nonlin} is devoted to  nonlinear problem.
We prove here that the problem generates  a dynamical system (see Theorem~\ref{th:wp}) which,
under some additional conditions, possesses a compact finite-dimensional
global attractor (Theorem~\ref{th:attractor}).

\section{Preliminaries}\label{sec:pre}
Let $D$ be a sufficiently smooth domain  in $\R^d$ and  $H^s(D)$ be the Sobolev space of order $s\in \R$
on  $D$ which we define (see \cite{Triebel78}) as restriction (in the sense of distributions)
 of the
space $H^s(\R^d)$ (introduced via Fourier transform).
We  define the norm in  $H^s(D)$ by the relation
\[
\|u\|_{s,D}^2=\inf\left\{\|w\|_{s,\R^d}^2\, :\; w\in H^s(\R^d),~~ w=u ~~
\mbox{on}~~D
    \right\}.
\]
We also use the notation $\|\cdot \|_{D}=\|\cdot \|_{0,D}$
for the corresponding $L_2$ norm and, similarly, $(\cdot,\cdot)_D$ for the $L_2$
inner product.
We denote by $H^s_0(D)$ the closure of $C_0^\infty(D)$ in  $H^s(D)$
(\wrt  $\|\cdot \|_{s,D}$) and introduce the spaces
\[
H^s_*(D):=\left\{f\big|_D\, :\;  f\in H^s(\R^d), \;
{\rm supp }\, f\subset \overline{D}\right\},\quad s\in \R,
\]
to describe
boundary traces on $\Om\subset\partial \cO$.
Since the extension by zero of elements from $H^s_*(D)$ give us
elements of $H^s(\R^d)$,
these spaces $H^s_*(D)$ can be treated not only as functional classes defined
on $D$ (and contained in $H^s(D)$) but also as  (closed) subspaces
of $ H^s(\R^d)$. We endow the classes $H^s_*(D)$ with the induced norms
 $\|f \|^*_{s,D}= \| f \|_{s,\R^d}$
for $f\in H^s_*(D)$. It is clear that
\[
\|f \|_{s,D}\le \|f \|^*_{s,D}, ~~ f\in H^s_*(D).
\]
It is known  (see \cite[Theorem 4.3.2/1]{Triebel78})
that $C_0^\infty(D)$ is dense in $H^s_*(D)$ and
\begin{align*}
& H^s_*(D)= H^s_0(D)~~\mbox{for}~~ -1/2<  s<\infty,~~ s-1/2\not\in
\{ 0,1,2,\ldots\}.
\end{align*}
The norms  $\|\cdot \|^*_{s,D}$  and $\|\cdot \|_{s,D}$ are equivalent
for these $s$.  Note that  in the notations of \cite{LiMa_1968}
the space $H^{m+1/2}_*(D)$ is the same as $H^{m+1/2}_{00}(D)$ for every
 $m= 0,1,2,\ldots$
Below we also use the factor-spaces
$H^s(D)/\R$  with the naturally induced norm.
\medskip\par
To describe fluid velocity fields
we first introduce  the class $\mathscr{C}_0(\cO)$ of
$C^\infty$ vector-valued solenoidal (i.e., divergence-free) functions
$ v=(v^1;v^2;v^3)$
on $\cO$ which vanish in a neighborhood  of $\pd\cO$ and also for $|x|$ large enough.
Then we denote by $\widetilde{X}$ the closure of $\sC_0(\cO)$ \wrt  the $L_2$-norm and
by $\widetilde{V}$ the closure  of $\sC_0(\cO)$  \wrt the $H^1$-norm. One
can see that
\[
\widetilde{X}=\left\{ v=(v^1;v^2;v^3)\in [L_2(\cO)]^3\, :\; {\rm div}\, v=0;\;
\gamma_n v\equiv (v,n)=0~\mbox{on}~ \pd\cO\right\}
\]
and
\begin{equation}\label{v-v-bar}
\widetilde{V}\subseteq \widetilde{V}^{\diamond}\equiv\left\{
v=(v^1;v^2;v^3)\in [H^1(\cO)]^3\, :\;  {\rm div}\, v=0; \;
v=0~\mbox{on}~ \pd\cO
  \right\}.
\end{equation}
For some details concerning this type of spaces see, e.g., \cite{lad-NSbook,temam-NS} and  \cite{galdi-book}.
\par
The following (geometry type) hypothesis plays an important role
in  our  further considerations.
\begin{assumption}[Domain Hypothesis]\label{as:dom}
{\rm We assume that
\begin{enumerate}
  \item[{\bf (i)}] there exists a smooth {\em bounded} domain  $\cO'\subseteq\cO$ such that
  $\overline{\Om}\subset\pd\cO'$;
  \item[{\bf (ii)}] we have the equality in \eqref{v-v-bar}, i.e., $\widetilde{V}=\widetilde{V}^{\diamond}$.
\end{enumerate}
}
\end{assumption}
 The sense of the first requirement in Assumption \ref{as:dom}
is obvious. As for the second one we refer to \cite{LadSol} for a discussion
of conditions on the domain which guarantee the equality
$\widetilde{V}=\widetilde{V}^{\diamond}$ (see also \cite[Sect.4.3]{galdi-book} and the references therein).
Here, as examples, we only note that this property holds in the  following cases:
(i) $\cO$ is a smooth domain  with the  compact boundary; (ii)
$\cO=\R^3_- =\{ x_3\le 0\}$; (iii) $\cO$ is given by \eqref{poisl-1}
with smooth bounded $\sB$ or with $\sB$ as in \eqref{lyier}
(infinitely long pipes and tubes of possibly varying cross section
are also admissible).
\par
We also need the Sobolev spaces consisting of functions with zero average
on the domain $\Om$, namely
we consider the subspace
\[
\widehat{L}_2(\Om)=\left\{u\in L_2(\Om): \int_\Om u(x') dx' =0 \right\}
\]
in $ L_2(\Om)$ and also the subspaces  $\widehat H^s(\Om)=H^s(\Om)\cap\widehat L_2(\Om)$
in $H^s(\Om)$ for $s>0$
with the standard $H^s(\Om)$-norm.
The notations   $\widehat H^s_*(\Om)$ and $\widehat H^s_0(\Om)$
have a similar meaning.
We denote by $\widehat{P}$   the projection
on  $\widehat{H}^2_0(\Om)$ in $H^2_0(\Om)$ which is orthogonal
with respect to the inner product $(\Delta\cdot, \Delta\cdot)_\Om$.
As it was already mentioned in \cite{ChuRyz2011} the subspace
$(I-\widehat{P})H^2_0(\Om)$
  consists of functions $u\in H^2_0(\Om)$ such that $\Delta^2u=const$ and
thus has dimension one.
\smallskip

In further considerations we need the following assertion concerning
extension of functions defined on $\Om$.

\begin{proposition}\label{pr:ext}
Let Assumption~\ref{as:dom}(i) be in force.
Then there exists a linear bounded operator $Ext\, :\, \widehat{L}_2(\Om)\mapsto \big[L_2(\cO)\big]^3$
such that
\[
{\rm div}\, Ext[\psi]=0~~\mbox{in}~\cO,~~ (Ext[\psi],n)\big|_{\Om}=\psi,~~(Ext[\psi],n)\big|_{S}=0,
\]
 and
 \begin{equation*}
  \| Ext[\psi] \|_{\big[H^{1/2-\delta}(\cO)\big]^3}\le C\|\psi\|_{\Om},~~~\forall \delta>0,~~
  \forall \psi\in\widehat{L}_2(\Om).
 \end{equation*}
Moreover,
\begin{itemize}
  \item if $\psi\in H^{s}_*(\Om)$ for some $0<s<1$, then $Ext[\psi]\in \big[H^{s+1/2}(\cO)\big]^3$
  with the estimate
 \begin{equation}\label{est-ext}
  \| Ext[\psi] \|_{\big[H^{s+1/2}(\cO)\big]^3}\le C\|\psi\|_{H^{s}_*(\Om)},
 \end{equation}
  and the relations
  $ Ext[\psi]\big|_{S}=(0;0;0)$ and $ Ext[\psi]\big|_{\Om}=(0;0;\psi)$ on the boundary of $\pd\cO$;
  \item  there exists    a smooth bounded subdomain  $\cO'$ in $\cO$ such that (i) $\Om\subset\pd\cO'$,
  (ii)~$ Ext[\psi]\big|_{\cO\setminus \cO'}=0$, and (iii) $Ext[\psi]\big|_{\cO'}\in \big[H^{2}(\cO')\big]^3$
provided $\psi\in H_0^{3/2+\delta}(\Om)$ for some $\delta>0$.
\end{itemize}
\end{proposition}
\begin{proof}
  On  a smooth bounded subdomain  $\cO'$ in $\cO$ such that  $\Om\subset\pd\cO'$
we consider the following Stokes problem:
\begin{align}
  -\nu\Delta v+\nabla p= 0, \quad
   \di v=0 \quad {\rm in}\quad \cO'; \nonumber
\\
 v=0 ~~ {\rm on}~\pd\cO'\setminus \Om;
\quad
v=(0;0;\psi) ~~{\rm on} ~ \Om,\label{stokes}
\end{align}
where  $\psi\in \Lto$ is given.
This type of boundary value problems  in \emph{bounded} domains
was studied by many authors  (see, e.g., \cite{lad-NSbook,temam-NS} and also the recent
monograph \cite{galdi-book}
and the references therein). To construct an extension operator we need the following properties of solutions
to \eqref{stokes}
(for some discussion and references concerning the assertion below we refer
to  \cite{ChuRyz2011}).

\begin{proposition}\label{pr:stokes}
Let  $\psi\in H^{s}_*(\Om)$ with
 $-1/2\le s\le  3/2$ and
 $\int_\Om\psi(x')dx'=0$.
Then
problem \eqref{stokes} has a unique solution
\[
\{v;p\} \in [H^{s+1/2}(\cO')]^3\times[ H^{s-1/2}(\cO')/\R]
\]
 such that
\begin{equation*}
\|v\|_{[H^{s+1/2}(\cO')]^3}+\|p\|_{H^{s-1/2}(\cO')/\R}
\le c_0 \|\psi\|_{H_*^{s}(\Om)}.
\end{equation*}
\end{proposition}
Now we can take  a solution $v$ to  \eqref{stokes} and  define $Ext[\psi]$ as the zero extension
of $v$ on the domain $\cO$. One can see that for this operator $Ext$ all statements
of Proposition~\ref{pr:ext} are in force.
\end{proof}
\begin{remark}\label{re:unbounded}
{\rm
We could not find in the literature an appropriate  statement of Proposition~\ref{pr:stokes}
for {\em unbounded} domains. On the other hand we do not know an
extension result in the class solenoidal functions with estimate \eqref{est-ext}
for some range of the parameter $s$.
This is why  we use this  way for a construction of the operator $Ext$.
We also note that in the case when $\cO$ is bounded we can take
$\cO'=\cO$. In this case $Ext$ is  a Green type operator
which maps $\psi$ into $v$ according to \eqref{stokes}.
Exactly this extension operator was used in \cite{ChuRyz2011}.
}
\end{remark}
Using the extension operator constructed above we introduce the spaces which we need
to describe the interaction between fluid and plate.
\par
Let Assumption~\ref{as:dom} be valid and
\[
\cM(\cO)=\left\{ v=v_0+Ext[\psi]\, :\; v_0\in \sC_0(\cO),~~ \psi\in \widehat{H}^2_0(\Om) \right\}.
\]
Then we denote by $X$ the closure of $\cM(\cO)$ \wrt  the $L_2$-norm and
by $V$ the closure  of $\cM(\cO)$  \wrt the $H^1$-norm. One
can see that
\[
X=\left\{ v=(v^1;v^2;v^3)\in [L_2(\cO)]^3\, :\; {\rm div}\, v=0;\;
\gamma_n v\equiv (v,n)=0~\mbox{on}~ S\right\}
\]
and
\begin{equation*}
V= V^{\diamond}\equiv\left\{
v=(v^1;v^2;v^3)\in [H^1(\cO)]^3\, \left| \begin{array}{l}
 {\rm div}\, v=0,\;
v=0~\mbox{on}~ S, \\ v^1=v^2=0~\mbox{on}~\Om \end{array} \right.
  \right\}.
\end{equation*}
We equip   $X$ with $L_2$-type norm $\|\cdot\|_\cO$
and denote by $(\cdot,\cdot)_\cO$ the corresponding inner product.
The space $V$ is endowed  with the standard $H^1$ norm.
\par
In conclusion of this section we mention that in the
the case of the Poiseuille flow in the tube
or between two planes  described above we deal with
a domain  satisfying the Friedrichs-P\'oincare
property\footnote{This property is valid in the case when the domain $\cO$ is bounded at least in one direction.}:
\begin{equation}\label{Fr-Poinc}
 \exists\, d_\cO>0: ~~  \int_\cO |v(x)|^2 dx\le d^2_\cO\int_\cO |\g v(x)|^2 dx,~~~\forall\, v\in H^1_0(\cO).
\end{equation}
By the localization argument
 one can show that
the inequality in \eqref{Fr-Poinc} implies a similar property
 for any $v\in\{g\in H^1(\cO):\, g|_{S}=0\}$ and thus
\begin{equation}\label{Fr-Poinc-V}
\exists\, c_\cO>0: ~~\|v\|_\cO\le c_\cO \|\nabla v \|_\cO,~~~ \forall\, v\in V,
\end{equation}
for the  Friedrichs-P\'oincare domains.
\par

\section{Linear problem}\label{sec:lin}
In this section we consider
 a linear version of \eqref{fl.1}--\eqref{Com-con}
which is  obtained
by replacing equation \eqref{pl_eq} with
its linear counterpart. Thus we deal with  the following
problem
\begin{align}\label{fl.1-lin}
&   v_t-\nu\Delta v+L_0v+\nabla p=G_f(t)~~\mbox{and}~~\di v=0~~{\rm in}~~
 \cO   \times(0,+\infty) \\[2mm]
& \label{fl.4-lin}
v=0 ~~ {\rm on}~S ~~\mbox{and}
~~
v\equiv(v^1;v^2;v^3)=(0;0;u_t) ~~{\rm on} ~~ \Om,
\\[2mm] &
u_{tt} + \De^2 u=G_{pl}(t)+p|_\Om  ~~{\rm on} ~~ \Om, \label{pl_eq-lin}
\\[2mm] &
u=\frac{\pd u}{\pd n} =0  ~~{\rm on} ~~ \pd\Om,
~~~
\int_\Om u_t(x',t) dx'=0 \quad \mbox{for all}~~ t\ge 0,
\label{plBC-lin}
\end{align}
which we
supply  with  initial data of the form
\begin{equation}
v(0)=v_0,\quad u(0)=u_0, \quad u_t(0)=u_1.  \label{IC-lin}
\end{equation}
Similarly to \cite{ChuRyz2011} we consider
 the following class  of test functions
\begin{equation*}
\cL_T=\left\{\phi \left|\begin{array}{l}
\phi\in L_2(0,T; \left[H^1(\cO)\right]^3),\; \phi_t\in L_2(0,T;  [L_2(\cO)]^3),  \\
{\rm div}\phi=0,\; \phi|_S=0,\; \phi|_\Om=(0;0;b),\; \phi(T)=0,  \\
b\in  L_2(0,T; \Hto),\; b_t\in  L_2(0,T; \Lto).
\end{array}\right.\right\}
\end{equation*}
and  introduce the following definition.
\begin{definition}\label{lin_de:solution}
{\rm
A pair of functions $(v(t);u(t))$ is said to be  a weak solution to
the problem in
\eqref{fl.1-lin}--\eqref{IC-lin}  on a time interval $[0,T]$ if
\begin{itemize}
    \item $v\in L_\infty(0,T;X)\bigcap L_2(0,T; V)$;
\item $u \in L_\infty(0,T;H^2_0(\Om)), \; u_t \in L_\infty(0,T; \Lto)$
and   $u(0)=u_0$;
    \item for every $\phi\in \cL_T$  the following equality holds:
        \begin{multline}
            -\!\int_0^T\!\!(v,\phi_t)_\cO dt  +\nu\!\int_0^T\!\!(\g v,\g\phi)_{\cO}  dt+\int_0^T\!\!(L_0v,\phi)_\cO dt \\ -\!\int_0^T\!\!(u_t,b_t)_\Om dt
           + \!\int_0^T\!\!(\De u, \De b)_\Om dt   \\
          =  \int_0^T(G_f, \phi)_{\cO} dt +\int_0^T(G_{pl},b)_\Om dt + (v_0, \phi(0))_{\cO} + (u_1, b(0))_\Om; \label{weak_sol_def}
        \end{multline}
    \item the compatibility condition  $v(t)|_\Om=(0;0;u_t(t))$
holds for almost all $t$.
   \end{itemize}
}
\end{definition}
 The same argument as
 in \cite{ChuRyz2011} shows
 that a weak solution $(v(t);u(t))$  satisfies the relation
       \begin{multline*}
           (v(t),\psi)_\cO
+(u_t(t),\beta)_\Om = (v_0, \psi)_{\cO} + (u_1, \beta)_\Om \\
-\int_0^t\Big[
 \nu(\g v,\g\psi)_{\cO} + (L_0v,\psi)_{\cO} \\ + (\De u, \De \beta)_\Om
          -  (G_f, \psi)_{\cO}  -(G_{pl},\beta)_\Om\Big] d\tau 
        \end{multline*}
for almost all $t\in [0,T]$ and for all
$\psi=(\psi^1;\psi^2;\psi^3)\in W$, where $\beta=\psi^3\big\vert_\Om$ and
\begin{equation}\label{space-W}
W=\left\{
\psi\in  V \left|  \;
  \psi|_\Om=(0;0;\beta),  \;
\beta\in \Hto \right. \right\}.
\end{equation}
It also follows from the compatibility condition and
the standard trace theorem that the plate velocity $u_t$ possesses an
additional spatial regularity, namely we have that
$u_t\in L_2(0,T; H^{1/2}_*(\Om))$.
\par
Below  as  phase spaces we  use
\begin{equation}\label{space-cH}
\cH=\left\{ (v_0;u_0;u_1)\in X\times H^2_0(\Om)\times\Lto :\; (v_0,n)\equiv
v_0^3 =u_1
~\mbox{on}~ \Om\right\}
\end{equation}
and
\begin{equation}\label{space-cH-hat}
\hch= \left\{ (v_0;u_0;u_1)\in \cH :\; u_0\in \wH^2_0(\Om)\right\}\subset\cH
\end{equation}
with the norm $\|(u_0;u_0;u_1)\|_\cH^2=\|v_0\|^2_{\cO}+\|\De u_0\|^2_{\Om}+\|u_1\|^2_{\Om}$.
\par
Our main result in this section is the following well-posedness
theorem concerning the linear problem.
\begin{theorem} \label{lin_WP}
Let Assumption~\ref{as:dom} be in force.
Assume that
\[
U_0=(v_0;u_0;u_1)\in \cH, ~~G_f(t)\in L_2(0,T; V'),~~G_{pl}(t)\in L_2(0,T; H^{-1/2}(\Om)).
\] Then
 for any interval $[0,T]$
there exists a unique weak solution $(v(t); u(t))$ to
\eqref{fl.1-lin}--\eqref{IC-lin}
 with the initial data $U_0$. This solution possesses the property
\begin{equation*}
U(t;U_0)\equiv U(t)\equiv (v(t); u(t); u_t(t))\in C(0,T; X\times H_0^2(\Om)\times \Lto),
\end{equation*}
and satisfies  the energy balance equality
\begin{multline}\label{lin_energy}
\cE_0(v(t), u(t), u_t(t))+\int_0^t\big[ \nu ||\g v||^2_\cO + (Av,  v)_\cO\big]d\tau \\ =\cE_0(v_0, u_0, u_1)
+\int_0^t(G_f,  v)_\cO d\tau +\int_0^t (G_{pl},u_\tau)_\Om d\tau
\end{multline}
 for every $t>0$, where the energy functional $\cE_0$ is defined
by the relation
\begin{equation*}
\cE_0(v(t), u(t), u_t(t))=\frac12\left(\|v(t)\|^2_\cO+ \|u_t(t)\|^2_\Om+\| \De u(t)\|_\Om^2\right).
\end{equation*}
\end{theorem}
If $G_f\equiv 0$ and $G_{pl}\equiv 0$, then  Theorem~\ref{lin_WP} implies that
the problem in \eqref{fl.1-lin}--\eqref{IC-lin} generates a strongly continuous semigroup.
In order to state our result on asymptotic stability of this semigroup
 we need additional assumptions.

\begin{assumption}[Stability Hypothesis]\label{as:stab}
{\rm
Assume that   one of the following conditions is valid:
\begin{itemize}
  \item either the matrix $A(x)$ in \eqref{oper-L0} is uniformly strictly
  positive, i.e.,
\[
\exists\, \si>0: (A(x)\xi,\xi)_{\R^3}\ge \si |\xi|_{\R^3}^2,~~ \forall \xi\in\R^3,\; x\in \overline{\cO};
\]
    \item or  the domain $\cO$ satisfies   the Friedrichs-P\'oincare
property \eqref{Fr-Poinc} and
\begin{equation}\label{FP-cond}
\exists\,\delta>0: (A(x)\xi,\xi)_{\R^3}\ge -
\left(\frac{\nu}{c^2_\cO} -\delta\right) |\xi|_{\R^3}^2,~~ \forall \xi\in\R^3,\; x\in \overline{\cO},
\end{equation}
 where $c_\cO$ is the constant from the Friedrichs-P\'oincare
inequality in \eqref{Fr-Poinc-V}.
\end{itemize}
}
\end{assumption}
Thus in the case of a general domain $\cO$ satisfying Assumption~\ref{as:dom}
to obtain a result on long-time dynamics we need to  assume the presence of some
additional damping mechanism (drag/friction terms). If the domain satisfies
the Friedrichs-P\'oincare property, then the result can be achieved
without any damping (e.g., we can take $A(x)\equiv 0$).
Moreover, we  note that the condition in \eqref{FP-cond}
is true when
$\sup_{x\in\cO}|A(x)|<\nu c^{-2}_\cO$,
where $|A(x)|$ is the operator
(Euclidian) norm in $\R^3$.
In the case of the  Poiseuille type flow (see \eqref{poisl-2})
this
means that
$|\g_{x_1,x_2}a|$ is small enough.
Since the profile $a$ can be written in the form $a=k\nu^{-1} a_*$,
where $a_*$ solves \eqref{poisl-1a} with $\nu=1$ and $k=1$,
the latter condition  is satisfied
when $ k\nu^{-2}\le c(\cB)$. Here $k$ is the  Poiseuille velocity
parameter and $c(\cB)$ is a constant depending on the cross-section $\cB$
the  tube $\cO$. In the case of the Oseen model
we have $a_0=(U;0;0)$ in \eqref{poisl-2} and thus there are  no
restrictions on the velocity $U$ of the underlining flow  for
Friedrichs-P\'oincare domains.

\begin{theorem}\label{th:stab}
In addition to the hypotheses of Theorem~\ref{lin_WP} we assume that
Assumption~\ref{as:stab} is in force.
Then
 there exist positive constants $M$ and $\ga$ such that
for every initial data $U_0=(v_0;u_0;u_1)$ from  $\hch$ we have
\begin{equation}\label{exp-st-nhm}
\|U(t)\|^2_\cH\le M e^{-\ga t} \|U_0\|^2_\cH+M \int_0^t
e^{-\ga(t-\tau)}\left[ \|G_f(\tau)\|_{V'}^2+
\|G_{pl}(\tau)\|^2_{-1/2,\Om}\right] d\tau.
\end{equation}
In particular,
if $G_f\equiv 0$ and $G_{pl}\equiv 0$, then the $C_0$-semigroup generated by \eqref{fl.1-lin}--\eqref{IC-lin}
is exponentially stable in $\hch$.
\par
In the  case of a general operator $L_0$
we need to add the term
\[
\left(M_1+M_2 \sup_{x\in\cO}|A(x)|\right) \int_0^t
e^{-\ga(t-\tau)} \|v(\tau)\|^2 d\tau
\]
in the right hand side of \eqref{exp-st-nhm}.
Here $|A(x)|$ denotes the operator
(Euclidian) norm in $\R^3$ and $M_1=0$ when
the domain $\cO$ satisfies   the Friedrichs-P\'oincare
property in \eqref{Fr-Poinc}.
\end{theorem}

\subsection*{Proof of Theorem~\ref{lin_WP}}
In the case when $L_0\equiv 0$ and $\cO$ is bounded this theorem was proved in \cite{ChuRyz2011}
(see also  \cite{OP1999} for a similar result).
We use the same idea as in  \cite{ChuRyz2011}.
 The main difficulty which we are faced
is that we loose several compactness  properties of the model
(e.g., we cannot use the basis of eigenfunctions of the Stokes operator).
\smallskip\par\noindent
{\em Step 1. Existence of an approximate solution.}
Let $\{\psi_i\}_{i\in \N}$ be an (orthonormal) basis in the space $\widetilde V$
consisting of the smooth finite in $\cO$ functions.
 Denote by $\{\xi_i\}_{i\in\N}$ the basis in $\Hto$
which consists of eigenfunctions of the following problem
\[
(\De \xi_i,\De w)_\Om=\kappa_i (\xi_i, w)_\Om,~~~\forall\, w\in \wH^2_0(\Om),
\]
with the eigenvalues $0<\kappa_1\le\kappa_2\le\ldots$ and $||\xi_i||_{\Om}=1$.
Let  $\phi_i=Ext[\xi_i]$, where the operator $Ext$
is defined in Proposition~\ref{pr:ext}. This proposition also yields
$\phi_i$ is $H^2$ in some vicinity of $\Om$ and thus as in \cite{ChuRyz2011}
 one can  conclude that $\pd_{x_3}\phi^3_i=0$ on $\Om$.
 \par
We define an approximate solution as a pair of functions
\begin{equation}
v_{n,m}(t)=\sum_{i=1}^m \alpha_i(t)\psi_i +\sum_{j=1}^n \dot{\beta}_j(t)\phi_j, \quad u_n(t)=\sum_{j=1}^n\beta_j(t)\xi_j + (I-\widehat{P})u_0, \label{approx_sol}
\end{equation}
satisfying the relations
   \begin{multline}
           (\dot v_{n,m}(t),\chi)_\cO
+(\ddot u_n(t), h)_\Om
+
 \nu(\g v_{n,m}(t),\g\chi)_{\cO} +(\De u_n(t), \De h)_\Om
\\
 = -(L_0v_{n,m},\chi)_\cO+
            (G_f(t), \chi)_{\cO}  +(G_{pl}(t),h)_\Om \label{app-sol1}
        \end{multline}
for  $t\in [0,T]$ and for
every $\chi$ and $h$ of the form
\begin{equation}\label{hi-h}
\chi=\sum_{k=1}^{m'} \chi_k\psi_k +Ext[h]~~\mbox{with} ~~
h=\sum_{k=1}^{n'} h_k\xi_k,
\end{equation}
where $m'\le m$ and $n'\le n$.
It is clear that $\chi\in W$ and $\chi\big|_\Om=(0;0;h)$.
The system in \eqref{app-sol1}  is endowed
with the initial data
\[
v_{v,m}(0)=\Pi_m(v_0-Ext[u_1])+Ext[P_nu_1],
\]
\[
 u_n(0)=P_n\widehat{P}u_0 + (I-\widehat{P})u_0, \; \dot{u}_n(0)=P_n u_1,
\]
 where $\Pi_m$  is  the
 orthoprojector  on $Lin\{\psi_j : j=1,\ldots,m,\}$ in $\widetilde{X}$ and $P_n$
is orthoprojector on
$Lin\{\xi_i : i=1,\ldots,n\}$ in $\Lto$.
Since $Ext:\Lto\mapsto X$,
it is clear that
\begin{equation*}
(v_{v,m}(0);u_n(0);  \dot{u}_n(0))\to  (v_0;u_0;u_1)~~
\mbox{strongly in $\cH$ as $m,n\to\infty$.}
\end{equation*}
\par
As in \cite{ChuRyz2011} one can show that (\ref{app-sol1})
can be reduced to some ODE in $\R^{m+n}$ and with given initial data has a unique  solution
on any time interval $[0,T]$.
\par
It follows from (\ref{approx_sol})  that
\[
v_{n,m}(t)=\sum_{i=1}^m \alpha_i(t)\psi_i + Ext[\pd_t u_n(t)].
\]
 This implies
the  boundary compatibility condition:
\begin{equation}\label{nm-comp}
v_{n,m}(t)=(0;0;\pd_t u_n(t))~~ \mbox{on}~~ \Om.
\end{equation}
{\em Step 2. Energy relation and a priori estimate
for an approximate solution.}
 Taking $\chi=v_{n,m}$ and $h=\pd_t u_n(t)$ in \eqref{app-sol1} we obtain
 the following energy balance relation for approximate solutions
\begin{align} \label{approx_est}
&\cE_0(v_{n,m}(t), u_n(t), \pd_t u_n(t))
\\ &{}\qquad
+ \nu \int_0^t\int_\cO |\g v_{n,m}|^2 dx d\tau
+\int_0^t(Av_{n,m},v_{n,m} )_\cO d\tau \nonumber
 \\ &  =\cE_0(v_{n,m}(0), u_n(0), \pd_t u_n(0)) +\int_0^t(G_f,v_{n,m} )_\cO d\tau + \int_0^t(G_{pl},\pd_t u_n)_\Om d\tau.
\nonumber
\end{align}
We use here the structure of $L_0$
which  after simple calculations (see, e.g., Lemma 1.3 \cite[Ch.2]{temam-NS}) yields the equality
$(L_0v_{n,m},v_{n,m} )_\cO=(Av_{n,m},v_{n,m} )_\cO$.
The relation in \eqref{approx_est} and Gronwall's lemma
implies the following a priori estimate
\begin{multline}
\sup_{t\in [0,T]}\left\{ \|v_{n,m}(t)\|^2_\cO+ \|\Delta u_n(t)\|^2_\Om +
 \|\pd_t u_n(t))\|^2_\Om\right\}
 \\
  +  \int_0^T \left(\|\g v_{n,m}\|_\cO^2 +  \|v_{n,m}\|^2_\cO\right) d\tau\le C_T. \label{approx_est1}
\end{multline}
By the trace theorem from (\ref{nm-comp}) and \eqref{approx_est1} we also have that
\begin{equation} \label{h12_est}
\int_0^T \|\pd_t u_n(\tau))\|^2_{H_*^{1/2}(\Om)} d\tau =
\int_0^T\| v_{n,m}(\tau)\|_{1/2,\pd\cO}^2 d\tau\le C_T.
\end{equation}
{\em Step 3. Limit transition.}
By \eqref{approx_est1} the sequence $\{(v_{n,m}; u_n; \pd_t u_n)\}$
contains a subsequence  such that
\begin{align}
&(v_{n,m}; u_n; \pd_t u_n) \rightharpoonup (v; u; \pd_t u) \quad \ast\mbox{-weakly in } L_\infty(0,T;\cH);\label{uv-conv} \\
&u_n \rightarrow u \quad \mbox{strongly in } C(0,T; H^{2-\e}_0(\Om))
,~~ \forall\, \eps>0;
\label{u-strong} \\
&v_{n,m} \rightharpoonup v \quad \mbox{weakly in } L_2(0,T;V).   \label{v_conv}
\end{align}
To obtain (\ref{u-strong}) we use the Aubin-Dubinsky  theorem
(see, e.g., \cite[Corollary~4]{sim}).
By (\ref{h12_est}) we can also suppose that
\begin{align}
&\pd_t u_n \rightharpoonup \pd_t u \quad \mbox{weakly in } L_2(0,T; H^{1/2}_*(\Om));\label{ut-conv}\\
&v_{n,m} \rightharpoonup v \quad \mbox{weakly in } L_2(0,T; H^{1/2}(\partial\cO)).
   \label{v_conv-b}
\end{align}
Applying the same argument as in \cite{ChuRyz2011}
and using relations \eqref{uv-conv}--\eqref{v_conv-b}
we conclude the proof of the existence of weak solutions which satisfy the
corresponding energy balance {\em inequality}. At this point we use Assumption~\ref{as:dom}(ii)
to approximate elements from $W$ by elements of the form \eqref{hi-h}.
We need this to establish \eqref{weak_sol_def} for $\phi\in\cL_T$.
\par
{\it Step 4. Uniqueness.}
We first consider the case when $L_0\equiv 0$
and use   Lions' idea (see \cite{Lions_1969}),
 with the same test function as  \cite{ChuRyz2011} in the case of a bounded domain.
After establishing properties of solutions in this case we consider the term $L_0v$ as a perturbation.
\par
Let $U^j(t)=(v^j(t);u^j(t);u_t^j(t))$, $j=1,2$, be two different solutions to the problem in question with the same initial data and $L_0\equiv 0$. Then their difference  $U(t)=U^1(t)-U^2(t)=(v(t);u(t);u_t(t))$ satisfies the variational equality
\begin{equation*}
-\int_0^T(v,\phi_t)_\cO  + \nu\int_0^T  (\g v, \g \phi)_\cO -\int_0^T  (u_t,  b_t)_\Om + \\
\int_0^T (\Delta u, \Delta b)_\Om =0 
\end{equation*}
for all $\phi\in\cL_T$, $b=(\phi|_\Om)^3$.
Now for every $0<s<T$ we take
\begin{equation*} 
\phi(t)\equiv\phi^s(t)=\left\{
    \begin{aligned}
         &-\int_t^s d\tau \int_0^\tau d\zeta v(\zeta), && t<s, \\
         &0, && t\ge s,
    \end{aligned}
\right.
\end{equation*}
as a test function. The same calculation as in  \cite{ChuRyz2011} yields the uniqueness in the case $L_0\equiv 0$.
\par
{\it Step 5. Continuity with respect to $t$ and the energy equality.}
Using the Lions lemma
(see \cite[Lemma 8.1]{LiMa_1968})
by the same argument as in \cite{ChuRyz2011}
 we first prove  any weak solution $(v(t);u(t);u_t(t))$ is weakly continuous in
$X\times H^2_0(\Om) \times \Lto$.

\par
To prove the energy equality (in the case $L_0=0$), we follow the scheme
of~\cite[Ch.1]{Lions_1969},
see also  \cite[Ch.3]{LiMa_1968}, in the form presented in  \cite{ChuRyz2011}.
Thus  as in \cite{ChuRyz2011} we can conclude that the solution is
strongly continuous in $t$. Moreover, the energy relation in the case $L_0=0$ with $G_f=0$
and $G_{pl}=0$ implies  that the corresponding solutions generates strongly continuous semigroup.
\par

{\it Step 6. Case $L_0\neq 0$}.
Using the energy relation for the problem with $L_0=0$ and $G_f(t):=G_f(t)-L_0v(t)$
we can establish the uniqueness  of solutions  via the Gronwall's type argument
and also the smoothness properties in the general case.
This  completes the proof of Theorem~\ref{lin_WP}.

\subsection*{Proof of Theorem~\ref{th:stab}}

 To prove the
estimate   in (\ref{exp-st-nhm}), we construct a Lyapunov function
using an idea from \cite{Chu_2010} (see also \cite{ChuRyz2011}). Let
\begin{equation*}
V(v_0,u_0,u_1)=\cE_0(v_0,u_0,u_1)+\e\Psi(v_0,u_0,u_1),
\end{equation*}
where
$\Psi(v_0,u_0,u_1)=(u_0,u_1)_\Om +(v_0,Ext[u_0])_\cO$
and $\e>0$ is
a small parameter which will be chosen later.
We consider these functionals on approximate solutions $(v_{n,m};u_n)$
for which $\widehat{P}u_0=u_0$ and thus  $\widehat{P}u_n(t)=u_n(t)$ for all $t>0$. This allow us to substitute in (\ref{app-sol1}) $Ext[u_n]$ instead of
$\chi$ and obtain that
  \begin{align}
\frac{d}{dt}\Psi_{n,m}(t) &\equiv
\frac{d}{dt}\Psi(v_{n,m}(t),u_{n}(t),
\pd_t u_{n}(t))
\notag \\
& =\|\pd_t u_{n}\|^2_\Om + (v_{n,m}, Ext[\pd_tu_n])_{\cO} -(L_0v_{n,m},Ext[u_n])_\cO
\notag \\
&
-
 \nu(\g v_{n,m},\g Ext[u_n])_{\cO} -\|\De u_n\|^2_\Om
      \notag    \\ &
           + (G_f, Ext[u_n])_{\cO}  +(G_{pl},u_n)_\Om. \label{app-psi}
        \end{align}
By Proposition~\ref{pr:ext}, using the compatibility
condition in (\ref{nm-comp}) and the trace theorem  we have that
\[
| (v_{n,m}, Ext[\pd_t u_n])_{\cO}|\le C \|v_{n,m}\|_\cO \|\pd_tu_n\|_{\Om}
\le C \left[ \|\g v_{n,m}\|^2_\cO+\| v_{n,m}\|^2_\cO\right].
\]
Similarly, for every $\eta>0$ we have
\[
  | (\g v_{n,m},\g Ext[u_n])_{\cO}|\le \eta \|\De u_n\|^2_\Om+ C_\eta
   \left[ \|\g v_{n,m}\|^2_\cO+\| v_{n,m}\|^2_\cO\right]
\]
and
\[
| (G_f, Ext[u_n])_{\cO}  +(G_{pl},u_n)_\Om|\le  \eta \|\De u_n\|^2_\Om+C_\eta \left[\|G_f\|_{V'}^2 +\|G_{pl}\|^2_{-1/2,\Om}\right].
\]
It is also clear that
\[
  | (L_0v_{n,m}, Ext[u_n])_{\cO}|\le \eta \|\De u_n\|^2_\Om+ C_\eta
  \left[ \|\g v_{n,m}\|^2_\cO+\| v_{n,m}\|^2_\cO\right].
\]
Therefore it follows from (\ref{app-psi}) that
\begin{align*}
\frac{d}{dt}\Psi_{n,m}(t) \le & -\frac12 \|\De u_n\|^2_\Om+ C  \left[ \|\g v_{n,m}\|^2_\cO+\| v_{n,m}\|^2_\cO\right]
\\ &
+C\left[\|G_f\|_{V'}^2 +\|G_{pl}\|^2_{-1/2,\Om}\right].
\end{align*}
Using the energy relation in (\ref{approx_est}) we also have that
\begin{multline*}
\frac{d}{dt}\cE_0(v_{n,m}(t), u_n(t), \pd_t u_n(t)) \le  -(\nu-\eta)
\|\g v_{n,m}\|^2_{\cO} +\eta\| v_{n,m}\|^2_\cO
\\
+ C_\eta \left[\|G_f\|_{V'}^2 +\|G_{pl}\|^2_{-1/2,\Om}\right]
 -(Av_{n,m},v_{n,m})_\cO,~~\forall\, \eta>0.
\end{multline*}
One can see that  the function
$V_{n,m}(t)\equiv V(v_{n,m}(t),u_{n}(t),
\pd_t u_{n}(t))$  satisfies the relations
\[
a_0 \cE_0(v_{n,m}(t),u_{n}(t),
\pd_t u_{n}(t))\le
V_{n,m}(t)\le a_1 \cE_0(v_{n,m}(t),u_{n}(t),
\pd_t u_{n}(t))
\]
for sufficiently small $\eps>0$. Using the stability hypothesis in Assumption~\ref{as:stab}
we can choose $\eta>0$ and $\si>0$ such that
\[
(\nu-\eta)
\|\g v_{n,m}\|^2_{\cO} -\eta\| v_{n,m}\|^2_\cO
+(Av_{n,m},v_{n,m})_\cO\ge \si \left[ \|\g v_{n,m}\|^2_\cO+\| v_{n,m}\|^2_\cO\right].
 \]
Therefore
we have that
\[
\frac d{dt} V_{n,m}(t)+a_2 V_{n,m}(t)\le a_3  \left[\|G_f\|_{V'}^2 +\|G_{pl}\|^2_{-1/2,\Om}\right]
\]
with positive constants $a_i$. This implies relation (\ref{exp-st-nhm})
for approximate solutions. The limit transition yields  (\ref{exp-st-nhm})
for every weak solution.
\par
In the general case we can apply (\ref{exp-st-nhm}) with $L_0:= (a_0,\g)v+\mu v$ and $G_f:=G_f-Av+\mu v$,
where $\mu>0$ (in the case of the Friedrichs-P\'oincare domains we can take $\mu=0$).
This implies the desired conclusion and completes the proof of Theorem~\ref{th:stab}.

\section{Nonlinear problem}\label{sec:nonlin}
In this section we deal with problem (\ref{fl.1})--(\ref{Com-con})
with a nonlinear feedback force.
First we   impose the following hypotheses concerning the force $\cF(u)$ in the plate equation (\ref{pl_eq}).

\begin{assumption}\label{A:force}
\begin{itemize}
 \item
There exists $\e>0$ such that $\cF(u)$ is locally Lipschitz from $H^{2-\e}_0(\Om)$ into $H^{-1/2}(\Om)$\footnote{
We recall \cite{Triebel78} that  $ H^{-1/2}(\Om) =[H^{1/2}_*(\Om)]'
 \varsupsetneqq [H^{1/2}_0(\Om)]'$.}
 in the sense that
\begin{equation}\label{f-lip}
\| \cF(u_1)-\cF_2(u_2)\|_{-1/2,\Om}\le C_R \| u_1-u_2\|_{2-\eps,\Om}
\end{equation}
for any $u_i\in H^{2}_0(\Om)$ such that $\| u_i\|_{2,\Om}\le R$.
  \item
  There exists a $C^1$-functional  $\Pi(u)$  on $H^2_0(\Om)$
such that $\cF(u)=\Pi'(u)$, where $\Pi'$ denotes the
 Fr\'echet derivative of  $\Pi$.
   \item
The plate force potential $\Pi$ is bounded on bounded sets from $H^2_0(\Om)$
 and  there exist $\eta<1/2$
 and $C\ge 0$ such that
\begin{equation}\label{8.1.1c1}
 \eta \|\De u\|_\Om^2 +\Pi(u)+C \ge 0\;,\quad \forall\, u\in H^2_0(\Om).
\end{equation}
\end{itemize}
\end{assumption}
Examples of
  nonlinear feedback (elastic) forces $\cF(u)$ satisfying Assumption~\ref{A:force}
are described in \cite{ChuKol2} and \cite{ChuRyz2011}, see also \cite{CLW-jde13}. They represent different plate models and include
  Kirchhoff,  von Karman, and Berger models.

\subsection{Well-Possedness}\label{sec:wp}

\begin{definition}\label{de:solution}
{\rm
A pair of functions $(v(t);u(t))$ is said to be  a weak solution to \eqref{fl.1}--(\ref{Com-con})  on a time interval $[0,T]$ if
\begin{itemize}
    \item $v\in L_\infty(0,T;X)\bigcap L_2(0,T; V)$;
\item $u \in L_\infty(0,T;H^2_0(\Om)), \; u_t \in L_\infty(0,T; \Lto)$,  $u(0)=u_0$;
    \item the  equality
in (\ref{weak_sol_def}) holds with $G_{pl}(t):= -\cF(u(t))+  G_{pl}(t)$;
    \item the compatibility condition  $v(t)|_\Om=(0;0;u_t(t))$ holds for almost all $t$.
\end{itemize}
}
\end{definition}

\begin{theorem}\label{th:wp}
Let Assumptions~\ref{as:dom} and \ref{A:force} be in force.
Assume that $U_0=(v_0;u_0;u_1)\in \cH$, $G_f(t)\in L_2(0,T; V')$
and  $G_{pl}(t)\in L_2(0,T; H^{-1/2}(\Om))$. Then
 for any interval $[0,T]$
there exists a unique weak solution $(v(t); u(t))$ to
\eqref{fl.1}--\eqref{Com-con}
 with the initial data $U_0$. This solution possesses the property
\begin{equation}\label{cont-ws-n}
U(t)\equiv (v(t); u(t); u_t(t))\in C(0,T; \cH),
\end{equation}
where $\cH$ is given by (\ref{space-cH}),
and satisfies  the energy balance equality
\begin{multline}\label{energy}
\cE(v(t), u(t), u_t(t))+ \int_0^t\left[\nu ||\g v||^2_\cO +(Av,  v)_\cO\right] d\tau
\\
=\cE(v_0, u_0, u_1)
+\int_0^t(G_f,  v)_\cO d\tau +\int_0^t (G_{pl}, u_\tau)_\Om d\tau
\end{multline}
 for every $t>0$, where the energy functional $\cE$ is defined
by the relation
\begin{equation*}
\cE(v, u, u_t)=\frac12\left(\|v\|^2_\cO+ \|u_t\|^2_\Om+\| \De u\|_\Om^2\right)+ \Pi(u).
\end{equation*}
Moreover, there exists a constant $a_{R,T}>0$ such that
  for any couple of weak solutions
  $U(t)=(v(t); u(t); u_t(t))$ and $\hat U(t)=(\hat v(t); \hat u(t); \hat u_t(t))$
with the initial data  possessing the property  $\|U_0\|_\cH, \|\hat U_0\|_\cH\le R$
we  have
\begin{equation}\label{mild-dif}
\|U(t)-\hat U(t)\|^2_\cH+\int_0^t\|\nabla( v-\hat v)\|_\cO^2 d\tau \le a_{R,T}
\|U_0-\hat U_0\|^2_\cH, ~~t\in [0,T].
\end{equation}
The spatial average of $u(t)$  is preserved. In particular, if $U_0\in\hch$,
then $U(t)\in\hch$ for every $t>0$.
We recall that $\hch$ is defined by (\ref{space-cH-hat}).
\end{theorem}

\begin{proof}
The proof of the local existence of an approximate solution is almost the same, as in the linear case (see Theorem \ref{lin_WP}).
We use approximate solutions of the same structure
 which satisfy (\ref{app-sol1}) with $-\cF(u_n(t))+  G_{pl}(t)$
instead of $G_{pl}(t)$.
Then using  the standard argument we establish
the energy relation in (\ref{energy}) for these approximate solutions.
Now
the positivity type estimate in (\ref{8.1.1c1}) allow us
to obtain  the same a priori estimates as
in (\ref{approx_est1}) and (\ref{h12_est}). Therefore
we can prove the global existence of approximate solutions and
 establish  the existence
of a weak solution $U(t)=(v(t); u(t); u_t(t))$ by the same
argument as in the linear case.
To  make  limit transition in the nonlinear term we
use (\ref{f-lip}).
\par
Next we  consider the  pair $(v(t); u(t))$ as a solution
to the linear problem with $G_{pl}(t):= -\cF(u(t))+  G_{pl}(t)$.
This allow us to obtain (\ref{cont-ws-n}) and also derive
energy balance relation  (\ref{energy}) from (\ref{lin_energy})  using
the potential structure of the force $\cF$:  $\cF(u)=\Pi'(u)$.
\par
Since the difference of two weak solution  can be treated as
a solution to the linear problem with $G_f\equiv 0$
and  $G_{pl}(t):= \cF(\hat u(t)) -\cF( u(t))$,
 we can obtain (\ref{mild-dif}) from the energy
equality  (\ref{lin_energy}).  The uniqueness
follows from  (\ref{mild-dif}).
\par
Preservation of the spatial average of $u(t)$
follows from the same property for approximate solutions.
\end{proof}

We can derive from  Theorem~\ref{th:wp} the following
assertion.
\begin{corollary}\label{co:generation}
In addition to the hypotheses of  Theorem~\ref{th:wp} we assume that
 $G_f(t)\equiv G_0\in V'$ is independent of $t$
 and $G_{pl}(t)\equiv 0$.
Then   problem \eqref{fl.1}--\eqref{Com-con} generates dynamical systems
 $(S_t, \cH)$ and $(S_t, \widehat{\cH})$ with the evolution operator
defined  by the formula $S_t U_0=(v(t);u(t);u_t(t))$,
where $(v;u)$ is a weak  solution to \eqref{fl.1}--\eqref{Com-con}
 with the initial data $U_0=(v_0; u_0; u_1)$.
 \par
 If we assume in addition that $G_0=0$ and Assumption~\ref{as:stab}  holds,
 then these systems are gradient with the full energy $\cE(v_0, u_0, u_1)$
as a Lyapunov function. This means that
(a) $U\mapsto \cE(U)$ is continuous on $\cH$,
(b)
 $\cE(S_t U_0)$ is not increasing in $t$, and
(c)
if $\cE(S_t U_0)=\cE(U_0)$ for some $t>0$, then
 $U_0$ is a stationary point of $S_t$
(i.e., $S_t U_0=U_0$ for all $t\ge 0$).
Moreover,
the set $\cE_R=\{U_0: \cE(U_0)\le R\}$ is a bounded
closed forward invariant set for every $R>0$.
\end{corollary}
\begin{proof}
The argument is the same as in  \cite{ChuRyz2011}.
We  only note that
under Assumption~\ref{as:stab}  from
(\ref{energy}) (with $G_f=0$ and $G_{pl}=0$) we have  that every stationary
point $U_*$  for $S_t$  has the form $U_*=(0; u; 0)$, where  $u\in H^2_0(\Om)$.
\end{proof}

\subsection{Stationary solutions}
As above we assume that $G_{pl}\equiv 0$ and $G_f\equiv 0$ and Assumptions~\ref{as:dom} and \ref{as:stab} holds.
It follows from Definition~\ref{de:solution} that a stationary
(time-independent) solution  is a pair $(v;u)$
 from $\widetilde{V}\times H_0^2(\Om)$
satisfying the relation
       \begin{equation}
 \nu(\g v,\g\psi)_{\cO} +(L_0v,\psi)_{\cO} +(\De u, \De \beta)_\Om
            +(\cF(u),\beta)_\Om=0
\label{weak_sol_stat}
        \end{equation}
for any $\psi\in W$ with $\psi^3\big\vert_\Om=\beta$,
where $W$ is given by (\ref{space-W}). Taking $\psi=v$ we conclude that
$\nu\|\g v\|^2_{\cO} +(Av, v)_{\cO}=0$ and hence from Assumption~\ref{as:stab}  we have
 $v=0$.
Therefore  we obtain the following variational problem for $u\in H^2_0(\Om)$:
      \begin{equation}
 (\De u, \De \beta)_\Om
        +(\cF(u),\beta)_\Om=0,~~~\forall\, \beta\in\wH^2_0(\Om).
\label{plate_sol_stat}
        \end{equation}
As in \cite{ChuRyz2011} we can  show the existence of a family of solutions to (\ref{plate_sol_stat})
parameterized by  a real parameter.
In the case of the zero average of $u$   we can fix this parameter and obtain  the following
assertion (see \cite{ChuRyz2011} for details).
\begin{proposition}[\cite{ChuRyz2011}]\label{pr:stat-sol}
In addition to Assumption~\ref{A:force} we
assume that
 there exist $\eta<1/2$  and  $c\ge 0$ such that
\begin{equation}\label{u-fu}
\eta \|\De u\|^2_{\Om}+ (u, \cF(u))_\Om \ge -c,~~~\forall\,
u\in H^2_0(\Om).
\end{equation}
Then  the set $\cN_0$ of solutions $u$
to  problem (\ref{plate_sol_stat}) with the property
$\int_\Om u dx=0$ is nonempty compact set in  $\wH^2_0(\Om)$.
This implies that the set $\cN$
of all stationary points of $S_t$ in the space $\hch$ is nonempty
compact set and has the form
\begin{equation}\label{N-set}
\cN=\left\{ (0;u;0): \; u\in  \wH^2_0(\Om)~~
\mbox{solves   (\ref{plate_sol_stat})}\right\}.
\end{equation}
\end{proposition}

\subsection{Asymptotical behavior} \label{sec:ab}

In this section we are interested in
global asymptotic behavior of the dynamical system  $(S_t, \hch)$.
Our main result states the existence of a compact global attractor
of finite fractal dimension.
 \par
We recall
(see, e.g., \cite{BabinVishik, Chueshov,Temam})
that
a \textit{global attractor}  of the  dynamical system  $(S_t, \hch)$
is defined as a bounded closed  set $\Ac\subset \hch$
which is  invariant ($S(t)\Ac=\Ac$ for all $t>0$) and  uniformly  attracts
all other bounded  sets:
$$
\lim_{t\to\infty} \sup\{{\rm dist}_\cH(S(t)y,\Ac):\ y\in B\} = 0
\quad\mbox{for any bounded  set $B$ in $\hch$.}
$$

\begin{theorem}\label{th:attractor}
Let Assumptions~\ref{as:dom}, \ref{as:stab},  and~\ref{A:force}  be in force.
Assume that $G_{pl}\equiv 0$, $G_f\equiv 0$ and (\ref{u-fu}) hold.
Then the   dynamical system  $(S_t, \widehat{\cH})$ possesses
a compact global attractor $\Ac$ of finite  fractal dimension\footnote{
For the definition and some  properties of
the {\em fractal dimension}, see, e.g., \cite{Chueshov} or \cite{Temam}.
}.
\par
 Moreover,
\begin{enumerate}
    \item[{\bf (1)}]
  Any trajectory
$\gamma=\{ (v(t);u(t); u_t(t)): t\in \R\}$
from the attractor $\Ac$
possesses the properties
\begin{equation}\label{u-smth}
(v_t;u_t;u_{tt})\in L_\infty (\R;  X\times \wH^2_0(\Om)\times \Lto)
\end{equation}
and there is $R>0$ such that
\begin{equation}\label{u-smth2}
\sup_{\ga \subset \Ac}\sup_{t\in\R}\left( \|v_t\|^2_{\cO}
+\|u_t\|^2_{2,\Om}+ \|u_{tt}\|^2_{\Om}\right)\le R^2.
\end{equation}
\item[{\bf (2)}]
The global attractor $\Ac$ consists of
full trajectories
$\{ (v(t);u(t); u_t(t)): t\in \R\}$
which are homoclinic to the set $\cN$, i.e.
\[
\lim_{t\to\pm\infty}\inf_{u_*\in\cN_0}\left( \|v(t)\|^2_{\cO}
+\|u- u_*\|^2_{2,\Om}+ \|u_{t}\|^2_{\Om}\right)=0,
\]
where  $\cN_0=\left\{ u\in  \wH^2_0(\Om)~~
\mbox{solves   (\ref{plate_sol_stat})}\right\}$. In addition we have
\begin{equation}\label{st-conv}
\lim_{t\to+\infty}\dist_{\widehat{\cH}} (S_ty, \cN)=0~~~\mbox{for any initial data $y\in \hch$.}
\end{equation}
\end{enumerate}
\end{theorem}
 We emphasize that Theorem~\ref{th:attractor}
deals with dynamics in the space $\hch$ (the case of the zero spatial average of the deflection).
  For a possible  approach to description of
the system long-time behavior in the space $\cH$
we refer  to \cite[Remark 4.9]{ChuRyz2011}.
\par
To obtain the result stated in Theorem~\ref{th:attractor}
it is sufficient to show
that the system is quasi-stable
 in the sense of Definition 7.9.2~\cite{cl-book}
 (see also Section 4.4 in \cite{cl-lect2012}).
For this   we can  use the stability properties of linear
problem (\ref{fl.1-lin})--(\ref{IC-lin}) established in
Theorem~\ref{th:stab}  and the same argument as in \cite{ChuRyz2011}
which yields   the following assertion.
\begin{lemma}[Quasi-stability] \label{pr:qst}
Let the hypotheses of Theorem~\ref{th:attractor} be in force and
$U^i(t)=(v^i(t);u^i(t);u^i_t(t))$, $i=1,2$,
 be two weak solutions with
initial data $U^i_0=(v^i_0;u^i_0;u^i_1)$ from  $\hch$
such that
 $\|U_0^i\|_\cH\le R$, $i=1,2$. Then
the difference
\[
Z(t)=U^1(t)-U^2(t)\equiv (v(t);u(t);u_t(t))
\]
satisfies the relation
\begin{equation}\label{q-stab}
\|Z(t)\|^2_\cH\le M_R e^{-\ga_* t} \|Z_0\|^2_\cH+M_R \int_0^t
e^{-\ga_*(t-\tau)}
\|u(\tau)\|^2_{\Om} d\tau
\end{equation}
for some positive constant $M_R$ and $\ga_*$.
\end{lemma}
\begin{proof} See \cite{ChuRyz2011} for some details.
\end{proof}
To complete
the proof  of Theorem~\ref{th:attractor} we note that by Proposition~7.9.4~\cite{cl-book}
 $(S_t, \widehat{\cH})$ is asymptotically smooth, i.e.,
 for any bounded set $B$ in $\widehat{\cH}$  such that $S_tB\subset B$
for $t>0$ there exists  a  compact set $K$ in the closure $\overline{B}$
of $B$, such that $S_tB$ converges uniformly to $K$.
 By Corollary~\ref{co:generation} the system is gradient.
 Proposition~\ref{pr:stat-sol}  yields
that
 the set $\cN$ of the stationary points (see \eqref{N-set}) is bounded in $\widehat{\cH}$.
Therefore  to prove  the existence of a  global attractor we can
 use  well-known  criteria   for gradient systems (see, e.g., \cite[Theorem 4.6]{Raugel}
 or Corollary 2.29 in \cite{cl-mem}).

The standard results on gradient systems with compact attractors
 (see, e.g., \cite{BabinVishik,Chueshov,Temam}) imply \eqref{st-conv}.
Since  $(S_t, \widehat{\cH})$
is quasi-stable,   the finiteness of fractal dimension ${\rm dim}_f\Ac$
follows from  Theorem~7.9.6~\cite{cl-book}.
 To obtain the result on regularity stated in (\ref{u-smth}) and
(\ref{u-smth2}) we apply  Theorem~7.9.8~\cite{cl-book}.


\end{document}